\documentclass{amsart}
\usepackage{amssymb}
\usepackage{amsthm}
\usepackage{epsfig}

\def\co{\colon\thinspace}
\def\R{\mathbb{R}}
\def\Z{\mathbb{Z}}

\def\N{\mathbb{N}}

\newcommand{\mc}[1]{\ensuremath{\mathcal{#1}}}

\renewcommand{\bar}[1]{\ensuremath{\overline{#1}}}
\newcommand{\llangle}{\langle\negthinspace\langle}
\newcommand{\rrangle}{\rangle\negthinspace\rangle}
\newcommand{\Parabs}{\{H_\lambda\}_{\lambda\in\Lambda}}
\newcommand{\Parabsprime}{\{H_i'\}_{i\in I}}

\newtheorem{theorem}{Theorem}
\newtheorem{lemma}[theorem]{Lemma}
\newtheorem{corollary}[theorem]{Corollary}
\newtheorem{proposition}[theorem]{Proposition}

\newtheorem{claim}{Claim}[theorem]
\newtheorem{subclaim}{Subclaim}[claim]

\theoremstyle{definition}
\newtheorem{remark}[theorem]{Remark}
\newtheorem{definition}[theorem]{Definition}

\begin{document}

\author[D. Groves]{Daniel Groves}
\address{California Institute of Technology}
\email{groves@caltech.edu}
\author[J. F. Manning]{Jason Fox Manning}
\address{University at Buffalo, SUNY}
\email{j399m@buffalo.edu}

\title[Fillings, finite generation and direct limits...]{Fillings, finite generation 
and direct limits of relatively hyperbolic groups}

\thanks{The first author was supported in part by NSF Grant DMS-0504251.
The second author was supported in part by an NSF Mathematical Sciences
Postdoctoral Research Fellowship.  Both authors thank the NSF for their support. 
We also thank Mark Sapir for asking us questions which lead to the results in
this paper, Guoliang Yu for help on the (Strong) Novikov Conjecture,
and the referee for several useful comments.}

\date{26 October, 2006}

\begin{abstract}
We examine the relationship between finitely and infinitely generated
relatively hyperbolic groups.
We observe that direct limits of relatively hyperbolic groups are
in fact direct limits of finitely generated relatively hyperbolic groups.
We use this (and known results) to prove
the Strong Novikov Conjecture for the groups constructed by Osin in \cite{Osin-small}.
\end{abstract}

\maketitle

\section{Introduction}
This paper is about the relationship between finitely generated and
infinitely generated relatively hyperbolic groups.  Most definitions
and characterizations
of relatively hyperbolic groups \cite{gromov,farb,bowditch,yaman,drutusapir} assume finite generation at
the outset.  One exception is that of Osin \cite{Osin-book}.  
In this definition, a
group, finitely presented relative to a system of subgroups, is
hyperbolic relative to that system if and only if the relative
isoperimetric function is linear. (See Section \ref{core} for definitions.)

Of interest to conjectures such as the Baum-Connes Conjecture and the
Novikov Conjecture are `exotic' groups, which may be used to test these
conjectures.  One way to build such exotic groups is by taking direct limits
of hyperbolic groups, a method brought to great prominence by
Olshanskii (see, for example, \cite{Olsh}).  However, hyperbolic groups
are known to satisfy the Baum-Connes Conjecture (this is due to
Mineyev and Yu \cite{MY}), and the rational injectivity of the Baum-Connes assembly
map is preserved by direct limits.  This implies the Novikov Conjecture
for groups built as direct limits of hyperbolic groups.

There are by now many results of the form:  Suppose that $\mc{U}$ is a 
property which implies the Novikov conjecture (or the Strong Novikov Conjecture,
or the Baum-Connes Conjecture), and suppose that $G$ is a (finitely
generated) group which
is hyperbolic relative to groups satisfying $\mc{U}$.  Then $G$ satisfies
$\mc{U}$.  (See, for example, \cite{Osin-FAD}, \cite{DG-ueRH}, \cite{Ozawa}).

It appears then that finitely generated relatively hyperbolic groups
are not a reasonable place to 
search for counterexamples to these conjectures, for likely one would
have to build the pathology into the parabolic subgroups, and relative
hyperbolicity would be irrelevant.  In \cite{Osin-small}, Osin
constructs groups with some remarkable properties as
direct limits of infinitely generated relatively hyperbolic groups.
One might think that these more flexible constructions might be of
some help in building counterexamples to the Strong Novikov Conjecture.
The main result of this paper (see Theorem \ref{Novikov}) shows that this
is probably not the
case: Behavior exotic enough to defy the Strong Novikov Conjecture
would have to be built into the finitely generated subgroups of the
parabolic subgroups.

An outline of this paper is as follows:  In Section \ref{core} we introduce the
notion of the {\em finitely generated core} of a relatively hyperbolic group
(this concept is implicit in \cite{Osin-book}).  In Section \ref{fillcore} we
explain how the main result in Osin's paper \cite{Osin-fill} about Dehn filling
in relatively hyperbolic groups follows from the version where the relatively
hyperbolic groups are assumed to be finitely generated.  (In \cite{GM} the
authors proved this finitely generated version under the additional assumption
that the group is torsion-free.)  In Section \ref{indlim} we record some
simple observations about direct limits of relatively hyperbolic groups.  In particular,
a direct limit of relatively hyperbolic groups is in fact the direct limit of finitely
generated relatively hyperbolic groups.  In Section \ref{s:Novikov} we use the
main result of Section \ref{indlim} (along with known results about uniform
embeddability and the Strong Novikov Conjecture) to deduce the main result
of this paper: that many of the remarkable groups constructed in \cite{Osin-small}
satisfy the Strong Novikov Conjecture (see Theorem \ref{Novikov}).

\section{The finitely generated core} \label{core}
Everything in this section is also contained in \cite{Osin-book}; only
the term {\em finitely generated core} is new.

Suppose that $G$ is a group with a family of subgroups
$\{H_\lambda\}_{\lambda\in\Lambda}$, so that $G$ has a \emph{finite
relative presentation}
\begin{equation}\label{relpres}
G = \langle S,\Parabs \mid \mc{R}\rangle,
\end{equation}
where $S$ is the finite relative generating set, and 
\[ \mc{R}\subset \mc{F}(S) \ast (\ast_{\lambda\in\Lambda}H_\lambda) \]
is a finite set of defining relations.  Here $\mc{F}(S)$ denotes the
free group on the letters $S$.  We briefly recall the meaning of
\eqref{relpres} from
\cite{Osin-book}.
Let $\mc{H} = \cup_{\lambda\in\Lambda}
H_\lambda$ be the disjoint union of the subgroups $\{H_\lambda\}$.
For each $\lambda$, let $\mc{R}_\lambda$ be the collection of words in
the alphabet $H_\lambda$ which represent the identity in $H_\lambda$.
The relative presentation 
\eqref{relpres} is shorthand for the (non-relative) presentation
\begin{equation}\label{nonrelpres}
G = \langle S\cup \mc{H}\mid \mc{R}\cup\mc{R}_{\Lambda}\rangle,
\end{equation}
where $\mc{R}_{\Lambda} = \cup_{\lambda\in\Lambda}\mc{R}_\lambda$.
If $D$ is a van Kampen diagram over the presentation
\eqref{nonrelpres}, then it has a \emph{relative area} which is the
number of $2$-cells in the diagram labelled by elements of $\mc{R}$.
Let $w$ be a word in the alphabet
$S\cup\mc{H}$ which represents the
identity in $G$.  The relative area of $w$ is the smallest relative
area of any van Kampen diagram over \eqref{nonrelpres} whose boundary
is labelled by $w$.
\begin{definition}
A group with finite relative presentation $G=\langle S,\Parabs \mid \mc{R}\rangle$ is \emph{hyperbolic relative to} the collection of subgroups
$\Parabs$ if the relative isoperimetric function is linear, in the
following sense:
There is a constant $C$ so that every word $w$ in the alphabet
$S\cup\mc{H}$ which represents the identity in $G$
can be filled by a van Kampen diagram of relative area
at most $Cn$.
\end{definition}

Given a relative presentation as in \eqref{relpres}, each $r\in \mc{R}$ can be
written in normal form as some product
\[ r = w_1 h_1\cdots w_n h_n\]
with each $w_i$ a word in the free group $\mc{F}(S)$ and each $h_i$ in
$H_{\lambda_i}$ for some $\lambda_i\in \Lambda$.  Let $\Omega$ be the
set of nontrivial elements of $\bigcup_\lambda H_\lambda$ which occur
in one of these normal forms.  Since $\Omega$ is finite, the set $I =
\{\lambda\in \Lambda\mid H_\lambda\cap \Omega\neq \emptyset\}$ is also
finite.  For each $i\in I$, $\Omega\cap H_i$ (finitely) generates some
$H_i'<H_i$.  Let $G'$ be the subgroup of $G$ generated by
$S\cup(\bigcup_{i\in I} H_i')$.
We call $G'$ the {\em finitely generated core of $G$ associated to the
presentation $ \langle S,\Parabs \mid \mc{R}\rangle$}.  (The finitely
generated core is highly non-unique; see Remark \ref{AddFiniteSet}.)
The next three lemmas are contained in the statement and proof of
Theorem 2.44 in \cite{Osin-book} (In Osin's statement, the finitely
generated core is called $Q$.). The second is straightforward; the
first and third can be proved using simple variations on the arguments in
Section \ref{fillcore}.
\begin{lemma}\label{l:corepres}
If $G'$ is the finitely generated core of $G$ coming from the finite
relative presentation $\langle S,\Parabs\mid \mc{R}\rangle$, then $G'$
has a finite relative presentation 
\begin{equation}\label{corepres}
 G' = \langle S,\Parabsprime\mid \mc{R}\rangle.
\end{equation}
(Note that
\[\mc{R}\subset \mc{F}(S) \ast (\ast_{i\in I}H_i') \subset \mc{F}(S) \ast
(\ast_{\lambda\in\Lambda}H_\lambda), \]
so the presentation in \eqref{corepres} is at least well defined.)
\end{lemma}
\begin{lemma}\label{l:freeproduct}
If $G$ is given by the finite relative presentation 
$\langle S,\Parabs \mid \mc{R}\rangle$, and $G' = \langle
S,\Parabsprime \mid\mc{R}\rangle$ is the corresponding finitely generated
core of $G$, as in the first paragraph, then 
\[G = G_0 \ast (\ast_{\lambda\in \Lambda\setminus I} H_\lambda)\]
where $G_0$ is the subgroup of $G$ generated by 
$S\cup (\bigcup_{i\in  I} H_i)$. 
\end{lemma}
\begin{lemma}\label{l:corerh}
$G= \langle S,\Parabs \mid \mc{R}\rangle$ is hyperbolic relative to
$\Parabs$ if and only if its finitely 
generated core $G'$ is hyperbolic relative to $\Parabsprime$.
\end{lemma}
\begin{remark} \label{AddFiniteSet}
If $G$ has a finite relative presentation
\[G = \langle S,\Parabs \mid \mc{R}\rangle\]
and $T\subseteq G$ is finite, then there is obviously a 
finite relative
presentation 
\[G = \langle S\cup T,\Parabs \mid \mc{R}'\rangle\]
for some (finite) set of relators $\mc{R}'$ containing $\mc{R}$.  Thus given
the finite set $T\subseteq G$, there is always a finitely generated
core containing $T$.  This observation will be used in Section \ref{indlim}.
\end{remark}

\section{Fillings and cores}\label{fillcore}
In \cite{Osin-fill}, Osin proved the following theorem:

\begin{theorem} \cite[Theorem 1.1]{Osin-fill} \label{infinite}
Suppose that $G$ is hyperbolic relative to the system of subgroups
$\{H_\lambda\}_{\lambda\in\Lambda}$. Then there is a finite set
$B\subset \cup_\lambda H_\lambda$ so that
if $\{K_\lambda\}_{\lambda\in\Lambda}$ is a collection of subgroups so
that for each $\lambda\in\Lambda$, we have
\begin{enumerate}
\item $K_\lambda\lhd H_\lambda$, and
\item $K_\lambda\cap B = \emptyset$,
\end{enumerate}
then 
\begin{enumerate}
\item The natural map $\phi_\lambda\co H_\lambda/K_\lambda\to
  G/\llangle \cup_\lambda K_\lambda\rrangle$ is an injection, and
\item $G/\llangle \cup_\lambda K_\lambda\rrangle$ is hyperbolic relative to
  $\{\phi_\lambda(H_\lambda/K_\lambda)\}_{\lambda\in\Lambda}$.
\end{enumerate}
\end{theorem}

In \cite{GM}, we proved the same theorem, with the 
additional assumptions that $G$ is torsion-free and that
the (finitely many) parabolic subgroups are finitely generated.

The torsion-free assumption in \cite{GM} is technical, and will not be
addressed here.  However, we remarked in \cite[Remark 1.5]{GM} that
if one can prove Theorem \ref{infinite} under the additional
assumptions that $\Lambda$ is finite and each $H_\lambda$ is finitely
generated, then the full statement follows.  This is the content
of  this section.

\begin{proposition} \label{finitetoinfinite}
If Theorem \ref{infinite} holds with the additional hypotheses that
$\Lambda$ is finite and that $H_\lambda$ is finitely generated for
each $\lambda\in\Lambda$, then Theorem \ref{infinite} holds in full
generality.
\end{proposition}

\begin{proof}
Let $G = \langle S,\Parabs\mid \mc{R}\rangle$ and $G'$ the finitely generated core of 
$G$ relative to $\langle S,\Parabs\mid \mc{R}\rangle$.  By Lemma \ref{l:corepres}, 
$\langle S,\Parabsprime\mid\mc{R}\rangle$ is a relative presentation for $G'$.
Lemma \ref{l:freeproduct} implies that $G$ splits as a free product
\[G = G_0 \ast (\ast_{\lambda\in \Lambda\setminus I} H_\lambda)\]
where $G_0$ is the subgroup of $G$ generated by 
$S\cup (\bigcup_{i\in  I} H_i)$.
Arbitrarily filling peripheral subgroups $H_\lambda$ for
$\lambda\notin I$ does not
affect the free product structure of $G$ (or relative hyperbolicity).
We may therefore assume that $\Lambda=I$, and so $G=G_0$.

Lemma \ref{l:corerh} implies that $G'$ is hyperbolic relative to $\{ H_i' \}_{i \in I}$,
so we may apply the finitely generated version of Theorem \ref{infinite} to $G'$.
Let $B$ be the finite subset of $\cup_\lambda H_\lambda'$ of
``forbidden'' elements for peripheral fillings of $G'$ coming from
Theorem \ref{infinite} in the finitely generated case.  We show that
$B$ also suffices as the set of forbidden elements for $G$.

For each $i\in I$, let $K_i\lhd H_i$ satisfy $K_i\cap B=\emptyset$.
The following is obvious:
\begin{claim}\label{claim:normal}
If $K_i' = K_i\cap H_i'$, then $K_i'\lhd H_i'$ and $K_i'\cap B=\emptyset$.
\end{claim}

Let $N$ be the normal closure of $\cup_i K_i$ in $G$, and let $N'$ be
the normal closure of $\cup_i K_i'$ in $G'$.  Write $\bar{H}_i$ for
the image of $H_i$ under the quotient
map $G\to G/N$, and write ${\bar{H}_i}'$ for the image of $H_i'$ 
under the quotient map $G'\to G'/N'$.
Theorem \ref{infinite} in the finitely generated case and Claim \ref{claim:normal} 
together imply
that $G'/N'$ is hyperbolic relative to $\{{\bar{H}_i}'\}_{i\in I}$.

We will show that $G'/N'$ is a finitely generated core of $G/N$;
Lemma \ref{l:corerh} and Claim \ref{injects} (below) then imply
the proposition.

We introduce some notation in this paragraph.  Let $\mc{H}$ be the
disjoint union of the $H_i$ for $i\in I$ and let $\mc{R}_I$
be the collection of words representing the identity in some $H_i$.
Finally, let $\mc{K}$ be the disjoint union of the $K_i$ for $i\in I$.
The quotient group $G/N$ has the presentation
\begin{equation}\label{quotientpres}
G/N=\langle S\cup\mc{H} \mid \mc{R}\cup\mc{R}_I\cup \mc{K} \rangle
\end{equation}
which can also be written as the relative presentation
\[
G/N = \langle S,\{H_i\}_{i\in I}\mid \mc{R}\cup \mc{K}\rangle.
\]
Let
$w\in\mc{F}(S)\ast(\ast_{i\in I} H_i)$ be a word which is trivial in
$G/N$.
There is a van Kampen diagram over \eqref{quotientpres}
for $w$ containing five possible kinds of two-cells:
\begin{enumerate}
\item \emph{$\mc{R}$-cells},
\item \emph{$H'$-cells}:  cells representing relations in $H_i'$ for some
  $i\in I$,
\item \emph{$K'$-cells}:  cells representing elements of
  $K_i'< G$ for some $i\in I$,
\item \emph{$H\setminus H'$-cells}: cells representing relations in $H_i$
  for some $i\in I$ which involve some elements of $H_i\setminus
  H_i'$, and
\item \emph{$K\setminus K'$-cells}: cells representing
  elements of $K_i\setminus K_i'$ for some $i\in I$.
\end{enumerate}
The first three kinds of $2$-cells will be called \emph{good}, the last
two \emph{bad}.  A \emph{bad patch} $P$ is a maximal union of bad
$2$-cells, subject to the condition that the interior of $P$ is connected.

Each edge of the van Kampen diagram is labelled by some element of
$S\cup(\bigcup_i H_i\setminus\{1\})$.  Edges of the van Kampen diagram
will be called \emph{good} if they are labelled by elements of
$S\cup(\bigcup_i H_i')$; otherwise they are \emph{bad}.  Every bad
$2$-cell has at least one bad edge in its boundary, whereas good
$2$-cells have no bad edges in their boundary.

The point of the good/bad notation is that if we can modify a van
Kampen diagram so that it contains only good $2$-cells, then it
follows that its boundary represents the trivial element of $G'/N'$.

\begin{claim}\label{primeinjects}
For each $i\in I$, the natural map from $H_i'/K_i'$ to $G/N$ is injective.
\end{claim}
\begin{proof}
Let $w\in H_i'$ be in the kernel of the map to $G/N$, and let $D$ be a
van Kampen diagram for $w$ as described above.  (In particular
$\partial D$ consists of a single good edge.)  If this van Kampen
diagram can be modified to contain only good $2$-cells, then $w$ is in
the kernel of the natural map from $H_i'$ to $G'/N'$, and thus by
the finitely generated version of Theorem \ref{infinite} we have $w\in K_i'$.

Let $P$ be any bad patch.
Note that all the $2$-cells in $P$ have boundary labels in
a single subgroup $H_j$, for some $j\in I$ which may be different from $i$.  
\begin{subclaim}\label{sc}
Each component of $\partial P$
represents an element of $H_j'$. 
\end{subclaim}
\begin{proof}
If not, then there
is some bad edge in $\partial P$.  This edge can only be
adjacent to another bad $2$-cell (in which case $P$ is not maximal) or
to the boundary of the van Kampen diagram itself.  Since the boundary
of the van Kampen diagram is labelled only by good edges, we derive a
contradiction.
\end{proof}

We now suppose $P$ is a bad patch which is \emph{innermost} in the
following sense:  No other bad patch is separated from the boundary of the van
Kampen diagram by the interior of $P$.  It is evident that if there are
any bad patches, then at least one is innermost.
Let $c_P$ be the outermost boundary of $P$.

We claim that the sub-diagram bounded by $c_P$ can be modified so that
every $2$-cell has all its boundary labels in $H_j$.  If $P$ is simply
connected, this is immediate.
Suppose then that the innermost bad patch $P$ is not simply connected.
Some component $c$ of $\partial P$ 
bounds a disk $D$ containing only good $2$-cells, since $P$ is
innermost.  Moreover, $c$ consists only of good edges.
Reading the labels of the edges of $c$ gives some word $w_c$ in the
alphabet $S\cup(\bigcup_i H_i')$, and the disk $D$ is itself a van
Kampen diagram for $w_c$ in $G'/N'$.  By Subclaim \ref{sc}, $w_c$
represents some element of $H_j'<G'$.  The disk $D$ is a demonstration
that it represents the trivial element of $G'/N'$.  Since $H_j'/N_j'$
injects into $G'/N'$ by assumption, it follows that $w_c$ represents
an element of $N'$.  The disk $D$ may therefore be replaced by a
single $K'$-cell.  Applying this argument in turn to each inner boundary
component of $P$, we fill $c_P$ entirely by $2$-cells with boundary
labels in $H_j$.

Let $w_P$ be the word in $S\cup(\bigcup_i H_i')$ given by the labels
of $c_P$.  By the previous paragraph, $w_P$ represents the trivial
element of $H_j/K_j$.  Moreover, by Subclaim \ref{sc}, it lies in
$H_j'$.  Since $H_j'/K_j'$ injects into $H_j/K_j$, $w_P$ represents an
element of $K_j'$.  We therefore can replace the bad patch $P$ (and
any disks attached to its interior) with a single $K'$-cell, thus
reducing the number of bad cells in the van Kampen diagram.  Iterating
this procedure, all bad cells can be removed, and Claim
\ref{primeinjects} follows.
\end{proof}
\begin{claim}\label{injects}
For each $i\in I$, the natural map from $H_i/K_i$ to $G/N$ is
injective.
\end{claim}
\begin{proof}
We now assume that $w\in H_i\setminus H_i'$ and suppose that $w$ lies
in $N$.  Again we can build a van Kampen diagram $D$ for $w$, this time
with boundary equal to a single bad edge.  We argue as in Claim \ref{primeinjects}:
First, there is a single outermost bad patch (since $\partial D$ is a a single bad edge),
and the boundary of each other bad patch is an element of some $H_j'$.  We now
reduce the number of bad patches until  $D$ contains 
a single bad patch and no $\mc{R}$-cells.  It follows that $w$ is already trivial
in $H_i/K_i$, and Claim \ref{injects} is proved.
\end{proof}

We now turn to the proof of Proposition \ref{finitetoinfinite}.
By Claim \ref{injects}, we can regard each $H_i/K_i$ as a subgroup of
$G/N$.  The subgroup $N$ of $G$ is the normal closure of the set
$\cup K_i$, so we obtain a finite relative presentation
\begin{equation}\label{quotient}
G/N = \langle S,\{H_i/K_i\}_{i\in I}\mid \bar{\mc{R}}\rangle,
\end{equation}
where $\bar{\mc{R}}$ is the image of $\mc{R}$ in 
$\mc{F}(S) \ast (\ast_{i\in I}H_i/K_i)$ under the obvious map from 
$\mc{F}(S) \ast (\ast_{i\in I}H_i)$.  It remains to observe that (regarding
$H_i'/K_i'$ as a subgroup of $H_i/K_i$) the group $G'/N'$ is the
finitely generated core coming from the presentation \eqref{quotient},
\[G'/N' = \langle S,\{H_i'/K_i'\}_{i\in I}\mid \bar{\mc{R}}\rangle.\]
This completes the proof of Proposition \ref{finitetoinfinite}.
\end{proof}

\section{Direct limits of relatively hyperbolic
  groups}\label{indlim}

Direct limits of groups are a particularly good way of building finitely
generated groups with interesting properties.  This idea was developed
by Olshanskii \cite{Olsh} with direct limits of hyperbolic groups, and 
recently by Osin \cite{Osin-small} for direct limits of relatively hyperbolic
groups.  We are interested in the relatively hyperbolic construction.

\begin{definition}
Suppose that $\{ G_i \}_{i \in \N}$ is a sequence of groups, and
$\{ \phi_i : G_i \to G_{i+1} \}$ is a sequence of homomorphisms. 

Let $X = \prod_{i \in \N} G_i$ be the Cartesian product of the $G_i$.
Define a subset $\Lambda \subset X$ as follows:
\[	\Lambda = \{ (g_i) \mid \exists J \ \forall j \ge J \ \  g_{j+1} = \phi_j(g_j) \}	.	\]
Put an equivalence relation `$\sim$' on $\Lambda$ so that $(g_i) = (h_i)$ if there
is some $K$ so that for all $k \ge K$ we have $g_k = h_k$.

Let $L = \Lambda / \sim$.  The group operation on $X$ descends to $L$,
and $L$ is a group, called the {\em direct limit of $\{ (G_i, \phi_i) \}$}.
\end{definition}

\begin{remark}
Using the definition above, if the groups $G_i$ are countable then so is the
limit $L$.  If one were to allow more general directed systems then this
would no longer be the case.  However, in order to understand the constructions
from \cite{Osin-small}, the above definition is sufficient.
\end{remark}

We will need the following lemma.

\begin{lemma} \label{l:countablepara}
Suppose that $G$ is a countable relatively hyperbolic group.  Then the
collection of (nontrivial) parabolic subgroups of $G$ is countable.
\end{lemma}
\begin{proof}
Let $\Lambda$ be the index set for the parabolic subgroups 
$\{ H_\lambda \}_{\lambda \in \Lambda}$ of $G$.  We suppose that
each $H_\lambda$ is nontrivial.  We have
to show that $\Lambda$ is countable.

Let $G'$ be a finitely generated core of $G$, and suppose that
$G'$ is hyperbolic relative to $\{ H_i' \}_{i \in I}$.

By Lemma \ref{l:freeproduct} we have
\[	G = G_0 \ast (\ast_{\lambda \in \Lambda \setminus I} H_\lambda),	\]
where $G_0$ is hyperbolic relative to $\{ H_i \}_{i \in I}$.  A countable group
cannot contain uncountably many nontrivial free factors, so $\Lambda \setminus I$
is countable, which implies that $\Lambda$ is countable, as required.
\end{proof}

\begin{proposition}\label{prop:limits}
Suppose that $G$ is a group which is isomorphic to a direct limit of groups
$\{G_i\stackrel{\phi_i}{\longrightarrow} G_{i+1}\mid i\in \N\}$ so that each $G_i$ is
countable and hyperbolic relative to some collection of 
proper subgroups $\mc{P}_i$.  Then for every $i$, there exists a finite relative presentation $Pres_i$ of $G_i$ with respect to $\mc{P}_i$ and homomorphisms
$\phi_i : G_i' \to G_{i+1}'$ such that $G$ is isomorphic to the
direct limit of \emph{finitely generated} groups
$\{G_i'\stackrel{\phi_i}{\longrightarrow} G_{i+1}'\mid i\in \N\}$, where $G_i'$
is the finitely generated core of $G_i$ associated to $Pres_i$ for each $i$.
\end{proposition}
\begin{proof}
Suppose that $G$ is the direct limit of $\{ (G_i , \phi_i ) \}$ where
$\phi_i : G_i \to G_{i+1}$, and the $G_i$ are all relatively hyperbolic.

Let $\mathcal{P}_i = \{ P^i_1, P^i_2, \ldots \}$ be the parabolic subgroups
of $G_i$, and $X_i$ a finite relative generating set for $G_i$ with respect
to $\mc{P}_i$.

The group $G$ is countable, so let $\{ g_0, g_1, g_2 , \ldots \}$ be an
enumeration of its elements.  For each $i \ge 0$, let $j(i)$ be the least
number so that (i) $j(i) \ge j(i-1) + 1$; and (ii) the image of $G_{j(i)}$ in $G$ contains
$\{ g_0 , \ldots , g_i \}$.  Since we may pass to a subsequence without changing
the limit, for ease of notation we will suppose that $j(i) = i$.

Let $Y_i = \{ g_0^{i}, \ldots , g_i^{i} \}$ be a subset of $G_i$
so that the map from $G_i$ to $G$ sends $g_l^{i}$ to $g_l$.

We will define a collection of finitely generated subgroups $G_i'$ of
$G_i$, and homomorphisms $\phi_i' : G_i' \to G_{i+1}'$ so that (i) the direct
limit of $\{ (G_i', \phi_i') \}$ is $G$; and (ii) each $G_i'$ is relatively hyperbolic.
In fact, the map $\phi_i'$ will be the restriction of $\phi_i$ to $G_i'$, and we
will use the notation $\phi_i$ for this map also.

Define $G_0'$ to be the finitely generated core of $G_0$ with respect
to a finite relative presentation with 
the relative generating set $Z_0 = X_0 \cup Y_0$.

Suppose that, for $r < i$, we have defined $G_r'$ (with finite relative 
generating set $Z_r$, and finitely many finitely generated parabolics) 
and $\phi_{r-1}'$.  We define
$G_i'$ as follows:

Let $W_i = \phi_{i-1}(Z_{i-1}) \in G_i$, and let $Z_i = W_i \cup X_i \cup Y_i$.
The set $Z_i$ is a finite relative generating set for $G_i$, so there
is some finite relative presentation $Pres_i = \langle
Z_i,\mc{P}_i\mid \mc{R}_i\rangle$.
Define $G_i'$ to be the finitely generated core of $G_i$ associated to
$Pres_i$.

By Lemma \ref{l:corerh}, the group $G_i'$ is hyperbolic relative to its
finitely many nontrivial intersections with the
elements of $\mc{P}_i$, and these intersections are themselves
finitely generated.

We claim that the direct limit of $\{ (G_i', \phi_i ) \}$ is isomorphic to $G$.
Let $L$ be the direct limit of $\{ (G_i', \phi_i )\}$.  
Since each $G_i'$ is a subgroup of $G_i$, and the map $\phi_i \co G_i' \to G_{i+1}'$
is a restriction of the homomorphism $\phi_i \co G_i \to G_{i+1}$, there is an obvious
map $\pi : L \to G$.  We construct the inverse map $\pi^{-1} \co G \to L$ as follows:
suppose that $g \in G$.  Then there is some $i$ so that $g = g_i$.  Then the choice
of $j(i)$ (and the renumbering above) implies that for all $j \ge i$ there is an element $g_i^j \in G_j$
so that $g_i^j$ maps to $g$ under the canonical map from $G_j$ to $G$.  But then
$g_i^j \in G_j'$ and $\phi_j(g_i^j) = g_i^{j+1}$.  We set $\pi^{-1}(g)$ to be the sequence
$(1,\ldots , 1, g_i^i, g_i^{i+1}, \ldots )$, where the first $i-1$
terms of this sequence are the identity element.  It is clear that $\pi^{-1}$ is the inverse of $\pi$, and we have
proved the proposition.
\end{proof}

\section{Some examples of Osin and the (Strong) Novikov conjecture} \label{s:Novikov}

In \cite{Osin-small}, Osin gives constructions of groups satisfying
some remarkable properties.  We show that if the input to these 
constructions is a uniformly embeddable group then the Strong Novikov
Conjecture\footnote{by which
we mean that the Baum-Connes assembly map is 
injective (see \cite{BCH,Yu-ICM}).} holds for the output.  In particular were these constructions to yield a 
counterexample to the Strong Novikov Conjecture, then the
input group must already have been rather exotic.

\emph{Uniform embeddability} for groups was introduced by Gromov in
\cite{Gromov-AsymInv}.  However, it is now more common to use the following
more general notion:

\begin{definition}\label{d:ue}
Let $(X,d_1)$ and $(Y,d_2)$ be metric spaces.  A map $i : X \to Y$ is a
{\em uniform embedding} if there are unbounded increasing functions 
$\rho_1, \rho_2 : \R^+ \to R^+$ so that for all $x,x' \in X$, 
\[	\rho_1(d_1(x,x')) \le d_2(i(x),i(x')) \le \rho_2(d_1(x,x'))	.	\]
A countable group $G$ is called {\em uniformly embeddable} if there is
a uniform embedding of $G$ into a Hilbert space.
\end{definition}
In the terminology of Roe  \cite[Chapter
  11]{RoeBook}, a group is uniformly embeddable if it admits a coarse
embedding into Hilbert space.
The connection between uniform
embeddability and the Novikov Conjecture was established in \cite{Yu},
where the Novikov Conjecture was proved for a uniformly embeddable group whose 
classifying space has the homotopy type of a finite CW complex.  The
finiteness assumption is removed in \cite{STY}, where the Strong
Novikov Conjecture is established for all uniformly embeddable groups.

In the statement below, $\pi(K)$ is the set of (finite) orders of
elements of a group $K$.  A group $G$ is said to be \emph{verbally
  complete} if the equation $w(x_1,\ldots,x_n)=g$ has a solution in $G$,
for any $g\in G$, and for $w$ any (freely reduced, nonempty) word in any
number of free 
variables. (In particular, every element of such a group is a
commutator, has roots of all orders, and so on.)  

\begin{theorem}\label{Novikov}
Given any uniformly embeddable countable group $G$, there exist $2$-generated
groups $H_1$ and $H_2$ satisfying:
\begin{enumerate}
\item $G$ embeds in $H_1$ and in $H_2$.
\item\label{Nov} $H_1$ and $H_2$ satisfy the Strong Novikov Conjecture.
\item $\pi(H_1)=\pi(G)$ and any two elements of $H_1$ with the same
  order are conjugate.
\item $H_2$ is verbally complete; moreover, if $G$ is torsion-free,
  then so is $H_2$.
\end{enumerate}
\end{theorem}
Before giving the proof of Theorem \ref{Novikov}, 
we should remark that the construction (and
most of the theorem) is due to Denis Osin in \cite{Osin-small}; the
sole innovation here is that point \eqref{Nov} can be guaranteed.

If we start with an infinite torsion-free uniformly embeddable group (like $\mathbb Z$),
then Theorem \ref{Novikov} yields:

\begin{corollary}
There exists a $2$-generated, infinite, torsion-free group $G$ which has two 
conjugacy classes and satisfies the Strong Novikov Conjecture.
\end{corollary}
\begin{corollary}
There exists a $2$-generated, infinite, torsion-free group $G$ which
is verbally complete and satisfies the Strong Novikov Conjecture.
\end{corollary}

The following proposition is used in the proof of Theorem \ref{Novikov}.
\begin{proposition}\label{prop:fad}
If $G$ is a uniformly embeddable group then $G$ can
be embedded into countable groups $R_1$ and $R_2$ so that:
\begin{enumerate}
\item Finitely generated subgroups of $R_1$ and $R_2$ are uniformly embeddable.
\item\label{orders} $\pi(R_1)=\pi(G)$ and all elements of the same order are
  conjugate.
\item\label{verb} $R_2$ is verbally complete.
\end{enumerate}
\end{proposition}
\begin{proof}
For $R_1$ satisfying \eqref{orders}, we use a construction of Higman,
Neumann, and Neumann \cite{HNN} (see Lyndon and Schupp \cite[Theorem
  IV.3.3]{LyndonSchupp}). Let $G_0=G$. 
Suppose $G_{i-1}$ has been defined, and let $\{(\alpha_j,\beta_j)\mid
j\in\N\}$ be the set of pairs of elements in $G_{i-1}$ so that the
orders of $\alpha_j$ and $\beta_j$ are equal.
Define $G_i$ by the presentation: 
\[G_i=\langle G_{i-1},\{t_{i,j}\}_{j\in\N}\mid t_{i,j}^{-1}\alpha_j
t_{i,j}=\beta_j, j\in \N \rangle, \]
and define $R_1 = \cup_{i\in \N} G_i$ to be the direct limit of these
groups.

For $R_2$ satisfying \eqref{verb}, we refer to the construction in
\cite{Osin-small}.  We will use only the following facts:
\begin{itemize}
\item The group $R_2$ is also a union of subgroups $U_i$ for $i\in
  \Z_{\geq 0}$, and $U_0=G$.  
\item If $i\geq 1$, then $U_{i+1}$ is an amalgamated free product of
$U_i$ with infinitely many groups $\{F_j^i\}_{j\in \N}$, where each
$F_j^i$ is either a free group or a one-relator group with torsion, and
each amalgamating subgroup is cyclic.
\end{itemize}

It remains only to prove that finitely generated subgroups of $R_1$
and $R_2$ are uniformly embeddable.  We prove both simultaneously. 
Let $H$ be a finitely generated subgroup of $R_1$ or $R_2$.  If $H<
R_1$, then we will set $V_i=G_i$ for each $i$; if $H<R_2$, we set $V_i=U_i$.
The group $H$ is contained in $V_i$ for some $i\in \N$.  If $i=0$, we
are done, since subgroups of uniformly embeddable groups are clearly
uniformly embeddable.  We may suppose by
induction that finitely generated subgroups of $V_{i-1}$ are uniformly
embeddable.   Free groups and one relator groups with torsion
are word hyperbolic (\cite{Newman}; see \cite[Theorem
  IV.5.5]{LyndonSchupp}); their finitely generated subgroups therefore
have finite asymptotic dimension
(this is a result of Gromov; see \cite{Roe-hypFAD}), and are therefore
uniformly embeddable (see \cite[Chapter 11]{RoeBook}).
Thus $V_i$ is a graph of groups with cyclic edge groups and
uniformly embeddable vertex groups.  The group $H$ inherits a graph of groups decomposition
from $V_i$; the edge groups are again cyclic.
Since $H$ is finitely generated, the graph of groups decomposition of $H$
has a finite underlying graph.  Also, since $H$ is finitely generated and the
edge groups of the (finite) graph of groups are finitely generated, the
vertex groups of the graph of groups are finitely generated.  Thus by
induction the vertex groups are uniformly embeddable.  Now we may
apply a theorem of Dadarlat and Guentner \cite{DG-ue} to
to deduce that $H$ is uniformly embeddable.
\end{proof}
\begin{remark}
Note that the proposition implies that $R_1$ and $R_2$ are uniformly
embeddable, since countable locally uniformly embeddable groups are uniformly
embeddable by \cite{DG-ue}.

The hypothesis and conclusion of
uniform embeddability in Proposition \ref{prop:fad}
may be strengthened to finite asymptotic
dimension, by applying theorems of Osin \cite{Osin-FAD} and Bell and
Dranishnikov \cite{BellDranisnikovTrees} in place of the theorems of
Dadarlat and Guentner quoted above.
\end{remark}

\begin{proof}(Theorem \ref{Novikov})
We do not describe Osin's constructions here, but refer the reader
to \cite{Osin-small}, particularly to the overview in Section 2 of his paper. 

The groups $H_1$ and $H_2$ are built as direct limits of (infinitely
generated) relatively hyperbolic groups, with peripheral subgroup
$R_1$ or $R_2$, respectively.  At each stage the term in the direct
limit is hyperbolic relative to $R_1$ or $R_2$ from the above
proposition. The finitely generated cores of these terms are therefore
hyperbolic relative to finitely generated subgroups of $R_1$ or
$R_2$; these finitely generated subgroups are uniformly embeddable by
Proposition \ref{prop:fad}.  The cores are thus hyperbolic relative to
uniformly embeddable subgroups; by the main result
of \cite{DG-ueRH} they are themselves uniformly embeddable.  
Thus, if $\Gamma$ is such a finitely generated core, then by
\cite[Theorem 6.1]{STY}, the Baum-Connes 
assembly map with coefficients is injective (with any separable
$\Gamma$-$C^\ast$-algebra coefficients).  
In particular, the Baum-Connes
assembly map (with trivial coefficients) is injective for each such
$\Gamma$, i.e. the finitely generated core $\Gamma$ satisfies the
Strong Novikov Conjecture. 

We have exhibited $H_1$ and $H_2$ as direct limits of groups
satisfying the Strong Novikov Conjecture.  By \cite[Proposition 2.4]{Rosenberg}, the Strong Novikov Conjecture is stable under taking
direct limits, so $H_1$ and $H_2$ themselves satisfy the Strong
Novikov Conjecture.
\end{proof}

By the main result of \cite{DG-ue}, the class of countable groups which are
uniformly embeddable is closed under direct limits of groups (where
all of the maps are injective).  It is clear from Lemma \ref{l:corerh} and
Remark \ref{AddFiniteSet} that
if $G$ is a countable relatively hyperbolic group with parabolic subgroups
$\{ H_i \}_{i \in I}$ then $G$ is the direct limit of an increasing
collection of finitely generated cores of $G$, each of which is hyperbolic
relative to a collection of finitely generated subgroups of finitely many of
the $H_i$.  Therefore, we have the following generalization of the
main result of \cite{DG-ueRH}.

\begin{proposition}
Suppose that $G$ is a (countable but not necessarily finitely generated)
group which is hyperbolic relative to a (not necessarily finite) collection
of subgroups $\{ H_i \}_{i \in I}$, and suppose that each of the $H_i$ is
uniformly embeddable.  Then $G$ is uniformly embeddable.
\end{proposition}

\begin{remark}
We feel that it is worth remarking that although Proposition \ref{prop:limits}
proves that in theory Osin could have built his examples without using infinitely
generated relatively hyperbolic groups, it is very difficult to see how to do
this directly, and in any case infinitely generated relatively hyperbolic groups
certainly make the proof conceptually easier.
\end{remark}

\small
\def\cprime{$'$} \providecommand\url[1]{\texttt{#1}}

\end{document}